\documentclass[11pt]{article}
\usepackage{graphics}
\usepackage[demo]{graphicx}
\usepackage{epsfig}
\usepackage{pstricks}
\usepackage[normal]{subfigure}
\usepackage[latin1]{inputenc}
\usepackage[english,francais]{babel}
\usepackage{relsize,exscale}
\usepackage{makeidx}
\usepackage{enumitem}
\usepackage{amsfonts,amssymb,amsmath}
\usepackage{graphicx}
\usepackage{color}
\usepackage{multirow}
\usepackage{mathrsfs}
\usepackage[normalem]{ulem}
\usepackage{cancel}
\usepackage{bbm}
\newenvironment{prooff}{{\it Proof :}}{\hfill\rule{2mm}{2mm}\vskip3mm\par}
\newtheorem{theorem}{Theorem}[section]
\newtheorem{lemma}[theorem]{Lemma}
\newtheorem{proposition}[theorem]{Proposition}

\newtheorem{e-definition}[theorem]{Definition\rm}
\newtheorem{remark}{\it Remark\/}

\usepackage{color}
\definecolor{dred}{rgb}{0.92,0,0}
\definecolor{dgreen}{rgb}{0,0.92,0}
\definecolor{dblue}{rgb}{0,0,0.92}
\definecolor{dyellow}{rgb}{0.95,0.95,0}

\newcommand{\R}{\mathbb{R}}
\newcommand{\N}{\mathbb{N}}
\def\D{\displaystyle}
\newcommand{\hs}{\hspace{0.1cm}}

\newcommand{\sa}{\\ [0.2cm]}
\usepackage{multirow}
\graphicspath{
{./Figures/}
{Figures/}
{./}
}
\title{A new first order Taylor-like theorem \\ with an optimized reduced remainder}
\author{Jo\"el Chaskalovic \thanks{D'Alembert,
Sorbonne University, Paris, France, (Email: \emph{(corresp.)} jch1826@gmail.com)}
\qquad
Hessam Jamshidipour
\thanks{
Sorbonne University, Paris, France, (Email: hessam.jamshidipour@gmail.com)}
}
\date{}
\begin{document}
\maketitle
\selectlanguage{english}
\begin{abstract}
\noindent This paper is devoted to a new first order Taylor-like formula where the corresponding remainder is strongly reduced in comparison with the usual one which appears in the classical Taylor's formula.
To derive this new formula, we introduce a linear combination of the first derivative of the concerned function, which is computed at $n+1$ equally-spaced points between the two points where the function has to be evaluated.
We show that an optimal choice of the weights in the linear combination leads to minimizing the corresponding remainder.
Then, we analyze the Lagrange $P_1$- interpolation error estimate and also the trapezoidal quadrature error, in order to assess the gain of accuracy we obtain using this new Taylor-like formula.
\end{abstract}
\noindent {\em keywords}: Taylor's theorem, Lagrange interpolation, interpolation error, quadrature trapezoid and corrected trapezoid rule, quadrature error.
\section{Introduction}\label{intro}
%
\noindent Rolle's theorem, and therefore, Lagrange and Taylor's theorems, prevent one from precisely determining the error estimate of numerical methods applied
to partial differential equations. Basically, this stems from the existence of a non-unique unknown point which appears in the remainder of Taylor's expansion,
as an heritage of Rolle's theorem. \sa
This is the reason why, in the context of finite elements, only asymptotic behaviors are generally considered for the error estimates, which strongly depend on the interpolation error
(see for example \cite{ChaskaPDE} or \cite{ArXiv_JCH}). \sa
Owing to this lack of information, several heuristic approaches were considered, so as to investigate new possibilities, which rely on a probabilistic approach. 
Such new possibilities finally enable one to classify numerical methods in which the associated data are fixed, and not asymptotic (for a review, see \cite{ArXiv_JCH}-\cite{MMA2021}) .\sa
However, an unavoidable fact is that Taylor's formula introduces an unknown point. This leads to the inability for one to determine exactly the interpolation 
error and, consequently, the approximation error of a given numerical method. It is thus legitimate to ask if the corresponding errors are bounded by quantities which 
are as small as possible.\sa
We wish here to focus on the values of the numerical constants which appear in these estimations so as to minimize them as much as possible. \sa
For example, let us consider the two-dimensional case and the $P_1$-Lagrange interpolation error of a given $C^2$ function which is defined on a given triangle. \sa
One can show that the numerical constant which naturally appears in the corresponding interpolation error estimate \cite{AsCh2014} is equal to $1/2$, as an heritage of the 
remainder of the first order Taylor expansion. \sa
This is the reason why, in this paper, we propose a new first-order Taylor-like formula, in which we strongly modify the repartition of the numerical weights between the Taylor polynomial and the corresponding remainder.\sa
To this end, we introduce a sequence of $(n+1)$ equally spaced points and consider a linear combination of the first derivative at these points. We show that an optimal
choice of the coefficients in this linear combination leads to minimizing the corresponding remainder, which becomes $2n$ smaller than the classical one obtained by the standard Taylor formula.\sa
As a consequence, we show that the Lagrange $P_1$-interpolation error as well as the quadrature error of the trapezoidal rule are two times smaller than the usual ones
obtained using the standard Taylor formula, provided we restrict ourselves to the new Taylor-like formula when $n=1$, namely, with two points.\sa
The paper is organized as follows. In Section \ref{B}, we present the main result of this paper, which regards the new first order Taylor-like formula.
In Section \ref{C}, Subsection \ref{Interpolation_Error}, we show the consequences we derived for the approximation error devoted to interpolation and,
in Subsection \ref{Quadrature_Error}, to that devoted to numerical quadratures. Finally, in Section \ref{D}, we provide concluding remarks.
\section{The new first-order, Taylor-like theorem}\label{B}
\noindent Let's first remind the well-known first order Taylor formula \cite{Taylor}. \sa
Let $(a,b) \in \mathbb{R}^{2}$, $a < b$, and $f \in \mathcal{C}^{2}([a,b])$. Then, it exists $(m_{2}, M_{2}) \in \mathbb{R}^{2} $ such that:
\begin{equation}
\label{m2M2}
\forall x \in [a,b] :\, m_{2} \leqslant f''(x) \leqslant M_{2}\,,
\end{equation}
and we have :
\begin{equation}\label{First_Order_Taylor}
f(b) = f(a) + (b-a)f'(a) + (b-a)\epsilon_{a,1}(b),
\end{equation}
where :
\[ \lim_{b \to a} \epsilon_{a,1}(b) = 0, \]
and:
\begin{equation}\label{Epsilon_b_bounded}
\frac{(b-a)}{2}m_{2} \leqslant \D \epsilon_{a,1}(b) \leqslant \frac{(b-a)}{2}M_{2}.
\end{equation}
In order to derive the main result below, let's introduce the function $\phi$ defined by:
\[
\begin{array}{r c c c l}
    \phi & : & [0,1] & \longrightarrow & \mathbb{R} \\
         &   &     t & \longmapsto & f'(a + t(b-a)).
\end{array}
\]
Then, we remark that : $\phi(0) = f'(a)$ and $\phi(1) = f'(b)$.
Moreover, the remainder $\epsilon_{a,1}(b)$ in (\ref{First_Order_Taylor}) satisfies the following result:
\begin{proposition}
The function $\epsilon_{a,1}(b)$ in the remainder (\ref{First_Order_Taylor}) can be written as follows:
\begin{equation}\label{Epsilon1}
\epsilon_{a,1}(b) = \int_{0}^{1}{(1-t)\phi'(t)dt}
\end{equation}
\end{proposition}
\begin{prooff}
Taylor's formula with the remainder in the integral form gives, at the first order:
\begin{equation}\label{Taylor_Integral}
f(b) = f(a) + (b-a)f'(a) + \frac{1}{1!}\int_{a}^{b}(b-x)f''(x)dx,
\end{equation}
and, using the substitution $x = a + (b-a)t$ in the integral of (\ref{Taylor_Integral}), we obtain :
\[ f(b) = f(a) + (b-a)f'(a) + (b-a)\int_{0}^{1}{(1-t)(b-a)f''(a + (b-a)t)dt,}\]
where :
\[ \phi'(t) = (b-a)f''(a + (b-a)t).\]
Finally,
\[ \epsilon_{a,1}(b) = \int_{0}^{1}{(1-t)\phi'(t)dt.}\]
\end{prooff}
Let now be $n \in \mathbb{N}^{*}$. We define $\epsilon_{a,n+1}(b)$ by the formula below:
\begin{equation}\label{Def_epsilon_n+1}
f(b) = f(a) + (b-a)\left(\sum \limits_{k=0}^{n} \omega_{k}(n)f'\left(a + k\frac{(b-a)}{n}\right)\right) + (b-a)\epsilon_{a,n+1}(b),
\end{equation}
where the sequence of the real weights $\D\biggl(\!\omega_{k}(n)\!\biggr)_{k\in [0,n]}$ will be determined such that the corresponding remainder built on $\epsilon_{a,n+1}(b)$ 
will be as small as possible.\sa
In other words, we will prove the following result:
\begin{theorem}\label{theorem_1}
Let $f$ be a real mapping defined on $[a,b]$ which belongs to $\mathcal{C}^{2}([a,b])$, such that: $\forall x \in [a,b], -\infty < m_2 \leqslant f''(x) \leqslant M_2 < +\infty$. \sa
If the weights $\omega_k(n), (k=0,n),$ satisfy: $$\D \sum \limits_{k=0}^{n} \omega_{k}(n) = 1,$$
then we have:
\begin{equation}\label{Poids}
\omega_{0}(n) = \omega_{n}(n) = \frac{1}{2n} \hs \mbox{ and } \hs \omega_{k}(n) = \frac{1}{n}, \hs \forall 0 < k < n,
\end{equation}
and the following first order expansion holds:
\begin{equation}\label{Generalized_Taylor}
f(b) = f(a) + (b-a)\left(\frac{f'(b) + f'(a)}{2n} + \frac{1}{n}\sum \limits_{k=1}^{n-1} f'\left(a + k\frac{(b-a)}{n}\right)\right) + (b-a)\epsilon_{a,n+1}(b),
\end{equation}
where :
\begin{equation}\label{majoration_epsilon_n+1}
\D\lim_{b \to a} \epsilon_{a,n+1}(b) = 0 \hs \mbox{ and } \hs |\epsilon_{a,n+1}(b)| \leqslant \frac{(b-a)}{8n}(M_2-m_2).
\end{equation}
Moreover, this result is optimal, since the weights given by (\ref{Poids}) guarantee that the remainder $\epsilon_{a,n+1}(b)$ is minimum.
\end{theorem}
\begin{remark}
To compare the control of $\epsilon_{a,n+1}(b)$ given by (\ref{majoration_epsilon_n+1}) and those of $\epsilon_{a,1}(b)$ given by (\ref{Epsilon_b_bounded}), we remark that (\ref{majoration_epsilon_n+1}) implies that:
\begin{equation}\label{majoration_epsilon_n+1_V2}
\D |\epsilon_{a,n+1}(b)| \leqslant \frac{(b-a)}{4n}\max(|m_2|,|M_2|).
\end{equation}
Consequently, the remainder $\epsilon_{a,n+1}(b)$ is $2n$ smaller than $\epsilon_{a,1}(b)$.
\end{remark}
\begin{remark}
We also notice in Theorem \ref{theorem_1} that between the parentheses lies a Riemann sum where, if $n$ tends to infinity, we obtain the fundamental theorem of integral calculus. That is to say :
\[ f(b) = f(a) + (b-a)\int_{0}^{1}{f'(bt + (1-t)a)dt}.\]
\end{remark}
In order to prove the theorem \ref{theorem_1}, we will need the following lemma:
\begin{lemma} \label{formula_1}
Let $u$ be any continuous function on $\R$, and a sequence of real numbers $(a_{k})_{0 \leqslant k \leqslant n}  \in \mathbb{R}^{n+1}, (n\in\N^*)$. We thus have the formula below :
\begin{equation}\label{Lemme_An_Bn}
 \sum \limits_{k=0}^{n-1} \int_{k}^{n}{a_{k} u(t) dt} = \sum \limits_{k=0}^{n-1} \int_{k}^{k+1}{S_{k} u(t) dt},
\end{equation}
where :
\[ S_{k} = \sum \limits_{j=0}^{k} a_{j}. \]
\end{lemma}
\begin{prooff}
Let's set : $\D A_{n} = \sum \limits_{k=0}^{n-1} \int_{k}^{n}a_{k}\,u(t)\, dt$ and $\D B_{n} = \sum \limits_{k=0}^{n-1} \int_{k}^{k+1}\,S_{k} u(t)\, dt$, where $\D S_{k} = \sum \limits_{j=0}^{k} a_{j}$.\sa
We will prove, by induction on $n$, that $A_{n} = B_{n}$ for all $n \in \mathbb{N}^{*}$.\sa
If $n = 1$, we indeed have :
\[ A_{1} = \sum \limits_{k=0}^{0} \int_{k}^{1}{a_{k} u(t) dt} = \int_{0}^{1}{a_{0} u(t) dt},\]
and :
\[ B_{1} = \sum \limits_{k=0}^{0} \int_{k}^{k+1}{S_{k} u(t) dt} = \int_{0}^{1}{S_{0} u(t) dt} = \int_{0}^{1}{a_{0} u(t) dt}.\]
So : $A_{1} = B_{1}$.\sa
Let us now assume that $A_{n} = B_{n}$, and let us show that : $A_{n+1} = B_{n+1}$.\sa
We have:
\begin{eqnarray}
\hspace{-0.6cm} B_{n+1} \hspace{-0.2cm}& = \hspace{-0.2cm}& \sum \limits_{k=0}^{n} \int_{k}^{k+1}{S_{k} u(t) dt} = \sum \limits_{k=0}^{n-1} \int_{k}^{k+1}{S_{k} u(t) dt} + \int_{n}^{n+1}{S_{n} u(t) dt} = B_{n} + \int_{n}^{n+1}{S_{n} u(t) dt}, \\[0.2cm]
& = & A_{n} + \int_{n}^{n+1}{S_{n} u(t) dt} = \sum \limits_{k=0}^{n-1} \int_{k}^{n}{a_{k} u(t) dt} + \int_{n}^{n+1}{S_{n} u(t) dt}, \\[0.2cm]
& = & \sum \limits_{k=0}^{n} \int_{k}^{n}{a_{k} u(t) dt} + \int_{n}^{n+1}{\left(\sum \limits_{j=0}^{n} a_{j}\right) u(t) dt} = \sum \limits_{k=0}^{n} \int_{k}^{n}{a_{k} u(t) dt} + \sum \limits_{k=0}^{n} \int_{n}^{n+1}{a_{k} u(t) dt}, \\[0.2cm]
 & = & \sum \limits_{k=0}^{n} \left( \int_{k}^{n}{a_{k} u(t) dt} + \int_{n}^{n+1}{a_{k} u(t) dt} \right) = \sum \limits_{k=0}^{n} \int_{k}^{n+1}{a_{k} u(t) dt} = A_{n+1}.
\end{eqnarray}
Conclusion :
 \[\forall n \in \N^*: \sum \limits_{k=0}^{n-1} \int_{k}^{n}{a_{k} u(t) dt} = \sum \limits_{k=0}^{n-1} \int_{k}^{k+1}{S_{k} u(t) dt},\]
where we set :
\begin{equation}\label{Sk(k)}
\D S_{k} = \sum \limits_{j=0}^{k} a_{j}.
\end{equation}
\end{prooff}
Let us now prove Theorem \ref{theorem_1}.\sa
\begin{prooff} We have :
\[ \D\frac{f(b) - f(a)}{b - a} = \phi(0) + \epsilon_{a,1}(b) = \sum \limits_{k=0}^{n} \omega_{k}(n)\phi\biggl(\frac{k}{n}\biggr) + \epsilon_{a,n+1}(b)\]
which can be re-written as:
\begin{eqnarray}
\epsilon_{a,n+1}(b) & = &\phi(0) + \epsilon_{a,1}(b) - \sum \limits_{k=0}^{n} \omega_{k}(n)\phi\biggl(\frac{k}{n}\biggr), \nonumber \\[0.2cm]
& = & \phi(0) + \int_{0}^{1}{(1-t)\phi'(t)dt} - \sum \limits_{k=0}^{n} \omega_{k}(n)\phi\biggl(\frac{k}{n}\biggr),\nonumber\\[0.2cm]
& = & \phi(1) - \int_{0}^{1}{t\phi'(t)dt} - \sum \limits_{k=0}^{n} \omega_{k}(n)\phi\biggl(\frac{k}{n}\biggr), \nonumber\\[0.2cm]
& = &\phi(1) - \int_{0}^{1}{t\phi'(t)dt} + \sum \limits_{k=0}^{n} \omega_{k}(n)\left(\phi(1)-\phi\biggl(\frac{k}{n}\biggr)\right) - \sum \limits_{k=0}^{n} \omega_{k}(n)\phi(1), \nonumber \\[0.2cm]
& = &\left(1 - \sum \limits_{k=0}^{n} \omega_{k}(n) \right)\phi(1) - \int_{0}^{1}{t\phi'(t)dt} + \sum \limits_{k=0}^{n} \omega_{k}(n)\left(\phi(1)-\phi\biggl(\frac{k}{n}\biggr)\right).\label{Equation espilon_n+1_1}
\end{eqnarray}
Let us assume for the sake of simplicity (see Remark \ref{Relation_Poids} below) that:
\begin{equation}\label{poids_norme}
\D \sum \limits_{k=0}^{n} \omega_{k}(n) = 1.
\end{equation}
(\ref{Equation espilon_n+1_1}) thus becomes:
\begin{eqnarray}
\epsilon_{a,n+1}(b) & = & - \int_{0}^{1}{t\phi'(t)dt} + \sum \limits_{k=0}^{n} \omega_{k}(n)\left(\phi(1)-\phi\biggl(\frac{k}{n}\biggr)\right),\nonumber \\[0.2cm]
 & = & - \int_{0}^{1}{t\phi'(t)dt} + \sum \limits_{k=0}^{n} \omega_{k}(n)\int_{\frac{k}{n}}^{1}{\phi'(t)dt}, \nonumber \\[0.2cm]
 & = & - \sum \limits_{k=0}^{n-1} \int_{\frac{k}{n}}^{\frac{k+1}{n}}{t\phi'(t)dt} + \sum \limits_{k=0}^{n} \omega_{k}(n) \int_{\frac{k}{n}}^{1}{\phi'(t)dt},\nonumber \\[0.2cm]
 & = & - \sum \limits_{k=0}^{n-1} \int_{\frac{k}{n}}^{\frac{k+1}{n}}{t\phi'(t)dt} + \sum \limits_{k=0}^{n-1} \int_{\frac{k}{n}}^{1}{\omega_{k}(n) \phi'(t)dt}.\label{Equation espilon_n+1_2}
\end{eqnarray}
Let us now use Lemma \ref{formula_1} in (\ref{Equation espilon_n+1_2}) by setting in (\ref{Lemme_An_Bn}):
$$\D u(t) = \phi'\biggl(\frac{t}{n}\biggr), \hs a_{k} = \omega_{k}(n), \mbox{ and } \hs S_{k} = \sum \limits_{j = 0}^{k} \omega_{j}(n) = S_{k}(n).$$
We thus obtain, using (\ref{Lemme_An_Bn}):
\[ \sum \limits_{k=0}^{n-1} \int_{k}^{n}{\omega_{k}(n) \phi'\biggl(\frac{t}{n}\biggr) dt} = \sum \limits_{k=0}^{n-1} \int_{k}^{k+1}S_{k}(n) \phi'\biggl(\frac{t}{n}\biggr) dt, \]
which can be written, by a simple substitution, as:
\begin{equation}
\D \sum \limits_{k=0}^{n-1} \int_{\frac{k}{n}}^{1}{\omega_{k}(n) \phi'(t) dt} = \sum \limits_{k=0}^{n-1} \int_{\frac{k}{n}}^{\frac{k+1}{n}}{S_{k}(n) \phi'(t) dt}.
\end{equation}
Then, (\ref{Equation espilon_n+1_2}) gives:
\begin{equation}\label{F1}
\D\epsilon_{a,n+1}(b) = - \sum \limits_{k=0}^{n-1} \int_{\frac{k}{n}}^{\frac{k+1}{n}}{t\phi'(t)dt} + \sum \limits_{k=0}^{n-1} \int_{\frac{k}{n}}^{\frac{k+1}{n}}{S_{k}(n) \phi'(t) dt} = \sum \limits_{k=0}^{n-1} \int_{\frac{k}{n}}^{\frac{k+1}{n}}{(S_{k}(n) - t)\phi'(t)dt}.
\end{equation}
Moreover:
$$\forall x \in [0,1] : \phi'(x) = (b-a) f''(a + x(b-a)), \hs \mbox{ and, } \forall t \in [a,b]: m_2 \leqslant f''(t) \leqslant M_2.$$
Next, to derive a double inequality on $\epsilon_{a,n+1}(b)$, we split the last integral in (\ref{F1}) as follows:
\begin{equation}
\int_{\frac{k}{n}}^{\frac{k+1}{n}}(S_{k}(n) - t)\phi'(t) dt = \int_{\frac{k}{n}}^{S_{k}(n)}(S_{k}(n) - t)\phi'(t) dt + \int_{S_{k}(n)}^{\frac{k+1}{n}}(S_{k}(n) - t)\phi'(t) dt.
\end{equation}
Then, considering the constant sign of $(S_{k}(n)-t)$ on $\biggl[\D \frac{k}{n}, S_{k}(n)\biggr]$, and on $\biggl[\D S_{k}(n),\frac{k+1}{n}\biggr]$, we have:
\begin{equation}\label{F2}
(b-a)m_2\int_{\frac{k}{n}}^{S_{k}(n)}(S_{k}(n) - t) dt \leq \int_{\frac{k}{n}}^{S_{k}(n)}(S_{k}(n) - t)\phi'(t) dt \leq (b-a)M_2\int_{\frac{k}{n}}^{S_{k}(n)}(S_{k}(n) - t) dt,
\end{equation}
and,
\begin{equation}\label{F3}
(b-a)M_2\int_{S_{k}(n)}^{\frac{k+1}{n}}(S_{k}(n) - t) dt \leq \int_{S_{k}(n)}^{\frac{k+1}{n}}(S_{k}(n) - t)\phi'(t) dt \leq (b-a)m_2\int_{S_{k}(n)}^{\frac{k+1}{n}}(S_{k}(n) - t) dt.
\end{equation}
Thus, (\ref{F2}) and (\ref{F3}) enable us to get the following two inequalities:
\begin{equation}\label{F4}
\int_{\frac{k}{n}}^{\frac{k+1}{n}}(S_{k}(n)-t)\phi'(t) dt \leq (b-a)M_2\int_{\frac{k}{n}}^{S_{k}(n)}(S_{k}(n)-t) dt + (b-a)m_2\int_{S_{k}(n)}^{\frac{k+1}{n}}(S_{k}(n)-t) dt,
\end{equation}
and,
\begin{equation}\label{F5}
\int_{\frac{k}{n}}^{\frac{k+1}{n}}(S_{k}(n)-t)\phi'(t) dt \geq (b-a)m_2\int_{\frac{k}{n}}^{S_{k}(n)}(S_{k}(n)-t) dt + (b-a)M_2\int_{S_{k}(n)}^{\frac{k+1}{n}}(S_{k}(n)-t) dt.
\end{equation}
Since we also have the following results:
\begin{equation}
\D \int_{\frac{k}{n}}^{S_{k}(n)}(S_{k}(n)-t) dt = \frac{\lambda^2}{2n^2} \hs \mbox{ and } \int_{S_{k}(n)}^{\frac{k+1}{n}}(S_{k}(n)-t) dt = -\frac{(\lambda-1)^2}{2n^2},
\end{equation}
where we set:
\begin{equation}\label{lambda}
\lambda \equiv nS_{k}(n)-k,
\end{equation}
inequalities (\ref{F4}) and (\ref{F5}) lead to:
\begin{equation}\label{F6}
\frac{(b-a)}{2n^2}P_1(\lambda) \leq \int_{\frac{k}{n}}^{\frac{k+1}{n}}(S_{k}(n)-t)\phi'(t) dt \leq  \frac{(b-a)}{2n^2}P_2(\lambda),
\end{equation}
where we defined the two polynomials $P_1(\lambda)$ and $P_2(\lambda)$ by:
\begin{equation}
P_1(\lambda) \equiv m_2\lambda^2-(\lambda -1)^2M_2 \hs \mbox{ and } \hs P_2(\lambda) \equiv M_2\lambda^2-(\lambda -1)^2m_2.
\end{equation}
Keeping in mind that we want to minimize $\epsilon_{a,n+1}(b)$, let us determine the value of $\lambda$ such that the polynomial $P(\lambda) \equiv P_2(\lambda)-P_1(\lambda)$ is minimum.\sa
To this end, let us remark that $P(\lambda)=(M_2-m_2)(2\lambda^2-2\lambda+1)$ is minimum when $\D\lambda = \frac{1}{2}$.\sa
Then, for this value of $\lambda$, (\ref{F6}) becomes:
\begin{equation}\label{F7}
\D \frac{(b-a)}{8n^2}(m_2-M_2) \leq \int_{\frac{k}{n}}^{\frac{k+1}{n}}(S_{k}(n)-t)\phi'(t) dt \leq  \frac{(b-a)}{8n^2}(M_2-m_2),
\end{equation}
and finally, by summing on $k$ between $0$ and $n$, we have:
\begin{equation}\label{epsilon1_final}
\frac{(b-a)}{8n}(m_2-M_2) \leqslant \epsilon_{a,n+1}(b) \leqslant \frac{(b-a)}{8n}(M_2-m_2).
\end{equation}
Due to definitions (\ref{Sk(k)}) of $S_{k}(n)$ and (\ref{lambda}) of $\lambda$, on the one hand, and because, on the other hand, the weights $\omega_{k}(n), (k=0,n)$ 
satisfy (\ref{poids_norme}),  we have:
\begin{equation}\label{Skn}
\forall n \in \mathbb{N}^{*}, \forall k \in [0,n[ : S_{k}(n) = \sum \limits_{j = 0}^{k} \omega_{j}(n) = \frac{1}{2n} + \frac{k}{n},
\end{equation}
and the corresponding $\omega_{k}(n)$ weights are equal to:
$$\omega_{0}(n) = \omega_{n}(n) = \frac{1}{2n}, \mbox{ and } \hs\omega_{k}(n) = \frac{1}{n},\, 0 < k < n,$$
which completes the proof for Theorem \ref{theorem_1}.
\end{prooff}
\noindent As an example, let us write formula (\ref{Generalized_Taylor}) when $n=2$ ({\sl i.e.}, with three points). In this case, we have:
\begin{equation}\label{Formule_3_points}
\D f(b) = f(a) + (b-a)\left(\frac{f'(a) + 2f'\biggl(\D\frac{a+b}{2}\biggr)+ f'(b)}{4}\right) + (b-a)\epsilon_{a,3}(b),
\end{equation}
where :
\[ \frac{(b-a)}{16}(m_2-M_2) \leqslant \epsilon_{a,3}(b) \leqslant \frac{(b-a)}{16}(M_2-m_2).\]
\begin{remark}\label{Relation_Poids}
Condition (\ref{poids_norme}) on the weights $\omega_k(n), (k=0,n),$ in Theorem \ref{theorem_1} is a kind of closure condition, since it helps determine $w_n(n)$,
but it is not a restrictive one. \sa
Indeed, without the closure condition (\ref{poids_norme}), one would have to consider the following expression of $\epsilon_{a,n+1}(b)$ instead of (\ref{F1}):
\begin{equation}\label{F1.1}
\D\epsilon_{a,n+1}(b) = \sum \limits_{k=0}^{n-1} \int_{\frac{k}{n}}^{\frac{k+1}{n}}{(S_{k}(n) - t)\phi'(t)dt} + \left(1 - \sum \limits_{k=0}^{n} \omega_{k}(n) \right)\phi(1).
\end{equation}
Then, (\ref{epsilon1_final}) would be replaced by:
\begin{equation}\label{epsilon1_final_V2}
\frac{(b-a)}{8n}(m_2-M_2) - \frac{M_1}{2n} \leqslant \epsilon_{a,n+1}(b) \leqslant \frac{(b-a)}{8n}(M_2-m_2) - \frac{m_1}{2n},
\end{equation}
where we assume that $(m_1,M_1) \in\R^2$ are such that $\forall x \in [a,b], -\infty < m_1 \leqslant f'(x) \leqslant M_1 < +\infty$. \sa
Moreover, to obtain (\ref{epsilon1_final_V2}), we also used the fact that the weights $\omega_k(n), (k=0,n)$ may be found with the help of (\ref{Skn}) without using anymore the 
closure condition (\ref{poids_norme}). \sa
More precisely, in this case, one can find that the weights $\omega_k(n), (k=0,n),$ are equal to:
\begin{equation}\label{Poids_New}
\omega_{0}(n) = \frac{1}{2n} \hs \mbox{ and } \hs \omega_{k}(n) = \frac{1}{n}, \hs \forall\, 0 < k \leq n.
\end{equation}
Consequently, from (\ref{epsilon1_final_V2}), we obtain that the remainder $\epsilon_{a,n+1}(b)$ is $n$ times smaller that those of the first-order Taylor formula given 
by (\ref{First_Order_Taylor}) and (\ref{Epsilon_b_bounded}).\sa
Finally, by considering the closure condition (\ref{poids_norme}) and the corresponding weights $\omega_k(n), (k=0,n)$, we slightly improved the result of (\ref{epsilon1_final_V2}),
since the remainder given by (\ref{F1}) is $2n$ smaller than those those of the first Taylor formula.
\end{remark}
\section{Application to the approximation error}\label{C}
\noindent To give added value to Theorem \ref{theorem_1}, which was presented in the previous section, this paragraph is devoted to appreciating the resulting differences one 
can observe in two main applications which belong to the field of numerical analysis. The first one concerns the Lagrange polynomial interpolation, 
and the second one the numerical quadrature. In these two cases, we will evaluate the corresponding approximation error both with the help of the standard first order 
Taylor formula and using the generalized formula (\ref{Generalized_Taylor}) derived in Theorem \ref{theorem_1}.
\subsection{The interpolation error}\label{Interpolation_Error}
\noindent In this section we consider the first application of the generalized Taylor-like expansion (\ref{Generalized_Taylor}) when $n=1$. In this case, for any function $f$ which belongs to $C^2([a,b])$, formula (\ref{Generalized_Taylor}) can be written :
\begin{equation}\label{Taylor_2pts}
\D f(b) = f(a) + (b-a)\left(\frac{f'(a)+f'(b)}{2}\right) + (b-a)\epsilon_{a,2}(b),
\end{equation}
where $\epsilon_{a,2}(b)$ satisfies :
\begin{equation}\label{inegalites_epsilon2}
\frac{(b-a)}{8}(m_2-M_2)\leqslant \epsilon_{a,2}(b) \leqslant \frac{(b-a)}{8}(M_2-m_2).
\end{equation}
\noindent As a first application of formula (\ref{Taylor_2pts})-(\ref{inegalites_epsilon2}), we will consider the particular case of the $P_1$-Lagrange interpolation (see \cite{Ciarlet2} or \cite{RaTho82}),
which consists in interpolating a given function $f$ on $[a,b]$ by a polynomial $\Pi_{[a,b]}(f)$ of degree less than or equal to one. \sa
Then, the corresponding polynomial of interpolation $\Pi_{[a,b]}(f)$ is given by:
\begin{equation}\label{P1_Interpolate}
\forall x \in [a,b]: \Pi_{[a,b]}(f)(x) = \frac{x-b}{a-b} \hs f(a) + \frac{x-a}{b-a} \hs f(b).
\end{equation}
One can remark that, using (\ref{P1_Interpolate}), we have: $\Pi_{[a,b]}(f)(a) = f(a)$ and $\Pi_{[a,b]}(f)(b) = f(b)$.\sa
Our purpose now is to investigate the consequences of formula (\ref{Taylor_2pts}) when one uses it to evaluate the error of interpolation $e(.)$ defined by $$\forall x \in [a,b]: e(x) = \Pi_{[a,b]}(f)(x) - f(x),$$ and to compare it with the classical first order Taylor formula given by (\ref{First_Order_Taylor}).\sa
Standard results \cite{Crouzeix_Mignot} regarding the $P_1-$Lagrange interpolation error claim that for any function $f$ which belongs to $C^2[a,b]$, we have:
\begin{equation}\label{error_interpolation_litterature}
\D|e(x)| \leq \frac{(b-a)^2}{2}\sup_{a\leq x\leq b}|f''(x)|.
\end{equation}
This result is usually derived by considering the suitable function $g(t)$ defined on $[a,b]$ by:
\begin{equation}\label{Error_bound}
g(t) = f(t)-\Pi_{[a,b]}(f)(t) - \biggl(f(t)-\Pi_{[a,b]}(f)(t)\biggr)\frac{(t-a)(t-b)}{(x-a)(x-b)}, (x \in ]a,b[).
\end{equation}
Given that: $g(a)=g(b)=g(x)=0$, and by applying Rolle's theorem two times, one can deduce that it exists $\xi_x \in ]a,b[$ such that: $g''(\xi_x)=0$.\sa
Therefore, after some calculations, one obtains the following:
\begin{equation}\label{N1}
f(x)-\Pi_{[a,b]}(f)(x) = \frac{1}{2 }(x-a)(x-b)f''(\xi), (a < \xi_x < b),
\end{equation}
and (\ref{error_interpolation_litterature}) simply follows. \sa
Still, as one can see from (\ref{N1}), estimation (\ref{error_interpolation_litterature}) can be improved since:
\begin{equation}\label{Max}
\D \sup_{a\leq x\leq b}(x-a)(b-x)=\frac{(b-a)^2}{4}.
\end{equation}
Then, (\ref{N1}) leads to:
\begin{equation}\label{error_interpolation_litterature_V2}
\D|e(x)| \leq \frac{(b-a)^2}{8}\sup_{a\leq x\leq b}|f''(x)|,
\end{equation}
in the place of (\ref{error_interpolation_litterature}).\sa
However, to appreciate the difference between the classical Taylor formula and the new one in (\ref{Taylor_2pts}), we will now reformulate the proof of
(\ref{error_interpolation_litterature_V2}) by using the classical Taylor formula (\ref{First_Order_Taylor}). This is the purpose of the following lemma:
\begin{lemma}\label{error_interpolation_classique_Thm}
Let $f$ be a function which belongs to $C^2([a,b])$ satisfying (\ref{m2M2}), then the first order Taylor theorem leads to the following interpolation error estimate:
\begin{equation}\label{error_P1_m_M}
|e(x)| \leq \frac{(b-a)^2}{8}M,
\end{equation}
where: $M=\max\{|m_2|,|M_2|\}$.
\end{lemma}
\begin{prooff}
We begin by writing the Lagrange $P_1-$polynomial $\Pi_{[a,b]}(f)$ given by (\ref{P1_Interpolate}) by the help of classical first order Taylor formula (\ref{First_Order_Taylor}). \sa
Indeed, in (\ref{P1_Interpolate}), we substitute $f(a)$ and $f(b)$ by:
\begin{eqnarray}
f(a) & = & f(x) + (a-x)f'(x) + (a-x)\epsilon_{x,1}(a), (\forall x \in [a,b]),\\ [0.2cm]
f(b) & = & f(x) + (b-x)f'(x) + (b-x)\epsilon_{x,1}(b), (\forall x \in [a,b]),
\end{eqnarray}
where, by the help of (\ref{Epsilon_b_bounded}) and (\ref{m2M2}), $\epsilon_{x,1}(a)$ and $\epsilon_{x,1}(b)$ satisfy:
\begin{equation}\label{epsilon_x,1(a)_x,1(b)}
\D|\epsilon_{x,1}(a)| \leq \frac{(x-a)}{2}M \hs \mbox{ and } \hs |\epsilon_{x,1}(b)| \leq \frac{(b-x)}{2}M.
\end{equation}
Then, (\ref{P1_Interpolate}) gives:
\begin{equation}\label{Error_P1_Taylor_V0}
\D \Pi_{[a,b]}(f)(x) = f(x) + \frac{(x-a)(b-x)}{(b-a)}\left[\epsilon_{x,1}(b)-\epsilon_{x,1}(a)\right],
\end{equation}
and due to (\ref{epsilon_x,1(a)_x,1(b)}), we get:
\begin{equation}\label{Error_P1_Taylor_V1}
\D |\Pi_{[a,b]}(f)(x) - f(x)| \leq \frac{(x-a)(b-x)}{2}M,
\end{equation}
where we used the fact that: $\D|\epsilon_{x,1}(b)-\epsilon_{x,1}(a)|\leq \frac{(b-a)}{2}M$.\sa
Finally, due to (\ref{Max}), (\ref{Error_P1_Taylor_V1}) leads to (\ref{error_P1_m_M}).
\end{prooff}
Let us now derive the corresponding result when one uses the new first order Taylor-like formula (\ref{Taylor_2pts}) in the expression of the interpolation polynomial $\Pi_{[a,b]}(f)$ defined by (\ref{P1_Interpolate}). \sa
This is the purpose of the following lemma:
\begin{lemma}\label{New_error_interpolation_Thm}
Let $f\in C^2([a,b])$, then we have the following interpolation error estimate:
\begin{equation}\label{error_P1_m_M_New}
\forall x \in [a,b]: \biggl|f(x) - \biggl[\Pi_{[a,b]}(f)(x) - \biggl(\frac{f'(b)-f'(a)}{2(b-a)}\biggr)(b-x)(x-a)\biggr]\biggr| \leq \frac{(b-a)^2}{32}(M_2-m_2).
\end{equation}
\end{lemma}
\begin{prooff}
We begin by writing $f(a)$ and $f(b)$ by the help of (\ref{Taylor_2pts}):
\begin{eqnarray}
f(a) & = & \D f(x) + (a-x)\left[\frac{f'(x)+f'(a)}{2}\right] + (a-x)\epsilon_{x,2}(a), \\[0.2cm]
f(b) & = & \D f(x) + (b-x)\left[\frac{f'(x)+f'(b)}{2}\right] + (b-x)\epsilon_{x,2}(b),
\end{eqnarray}
where $\epsilon_{x,2}$ satisfies (\ref{inegalites_epsilon2}), with obvious changes of notations. Namely, we have:
\begin{equation}\label{epsilon_x,2(a)_x,2(b)}
\D|\epsilon_{x,2}(a)| \leq \frac{(x-a)}{8}(M_2-m_2) \hs \mbox{ and } \hs \D|\epsilon_{x,2}(b)| \leq \frac{(b-x)}{8}(M_2-m_2).
\end{equation}
Then, by substituting $f(a)$ and $f(b)$ in the interpolation polynomial given by (\ref{P1_Interpolate}), we have:
\begin{equation}\label{Error_P1_V00}
\D \Pi_{[a,b]}(f)(x) = f(x) + \biggl(\frac{f'(b)-f'(a)}{2(b-a)}\biggr)(b-x)(x-a) + \frac{(b-x)(x-a)}{(b-a)}\left[\frac{}{}\!\epsilon_{x,2}(b)-\epsilon_{x,2}(a)\right].
\end{equation}
Now, if we introduce the \emph{improved} interpolation polynomial $\Pi^*_{[a,b]}(f)$ defined by:
\begin{equation}\label{corrected_polynomial}
\forall x \in [a,b]: \Pi^*_{[a,b]}(f)(x) = \Pi_{[a,b]}(f)(x) - \biggl(\frac{f'(b)-f'(a)}{2(b-a)}\biggr)(b-x)(x-a),
\end{equation}
equation (\ref{Error_P1_V00}) becomes:
\begin{equation}\label{Error_P1_V00_1}
\D \Pi^*_{[a,b]}(f)(x) = f(x) + \frac{(b-x)(x-a)}{(b-a)}\left[\frac{}{}\!\epsilon_{x,2}(b)-\epsilon_{x,2}(a)\right].
\end{equation}
Thus, due to (\ref{epsilon_x,2(a)_x,2(b)}), we have: $\D|\epsilon_{x,2}(a)-\epsilon_{x,2}(b)|\leq \frac{(b-a)}{8}(M_2-m_2)$, and (\ref{Error_P1_V00}) by the help of (\ref{corrected_polynomial}) gives:
\begin{equation}
\D |\Pi^*_{[a,b]}(f)(x) - f(x)| \leq \frac{(b-x)(x-a)}{8}(M_2-m_2) \leq \frac{(b-a)^2}{32}(M_2-m_2),
\end{equation}
which completes the proof of this lemma.
\end{prooff}
Let us now formulate a couple of consequences of Lemma  \ref{error_interpolation_classique_Thm} and Lemma \ref{New_error_interpolation_Thm}:
\begin{enumerate}
\item If we consider the improved interpolation polynomial $\Pi^*_{[a,b]}(f)$ defined by (\ref{corrected_polynomial}), we obtain for the error estimate 
(\ref{error_P1_m_M_New}) an accuracy which is two times more precise than that we got in (\ref{error_P1_m_M}) by the classical Taylor formula. \sa
    In order to compare (\ref{error_P1_m_M}) and (\ref{error_P1_m_M_New}), we notice that (\ref{error_P1_m_M_New}) leads to:
    \begin{equation}\label{Max2}
    \forall x \in [a,b]: |e^*(x)| \leq \frac{(b-a)^2}{16}\max(|m_2|,|M_2|).
    \end{equation}
    Now, the overcost for this improvement is that $\Pi^*_{[a,b]}(f)$ is a polynomial of degree less than or equal to two which requires the computation of $f'(a)$ and $f'(b)$. However, the consequent gain clearly appears in the following application devoted to finite elements. \sa
    Indeed, due to C\'ea's lemma \cite{ChaskaPDE}, the approximation error is bounded by the interpolation error. Thus, if one wants to locally guarantee that 
    the upper bound of the interpolation error be not greater than a given threshold $\epsilon$, then: \sa
    if $h$ denotes the local mesh size defined by $h=b-a$, with the classical $P_1$-Lagrange interpolation, and by $h^*$ the corresponding
    one with the improved interpolation $\Pi^*_{[a,b]}(f)$, we have:
    \begin{equation}\label{Perf_Interpol_New}
    \D \frac{M}{8}h^2 \leq \epsilon \mbox{ and } \D \frac{M}{16}h^{*2} \leq \epsilon.
    \end{equation}
It follows that the difference of the maximum size between $h$ and $h^*$ is $\D 1/\sqrt{2}\simeq 0.707$ and, consequently, $h^*$ may be around $30$ percents greater that $h$. This economy in terms of the total
number of meshes would be even more significant if one considers the extension of this case to a three dimensional application.
\item We also notice that, if we now consider the particular class of $C^2-$functions $f$ defined on $\R$, $(b-a)$-periodic, then $f'(a)=f'(b)$ and, consequently, the interpolation error $e^*(x)$ is equal to $e(x)$, and (\ref{error_P1_m_M_New}) becomes:
\begin{equation}\label{Error_P1_2pts_V3}
\forall x \in [a,b]: |e^*(x)| = |e(x)| \leq \frac{(b-a)^2}{32}(M_2-m_2) \leq \frac{(b-a)^2}{16}M.
\end{equation}
In other words, for this class of periodic functions, due to the new first order Taylor-like formula (\ref{Taylor_2pts}), the interpolation error $e(x)$ provided by 
(\ref{Error_P1_2pts_V3}) is bound by a quantity which is two times smaller that those we got in (\ref{error_P1_m_M}) using the classical Taylor formula. \sa
We highlight that in this case, there is no overcost anymore to obtain this more accurate result, since it concerns the standard interpolation error associated with the 
standard Lagrange $P_1$-polynomial.
\item Finally, since the improved polynomial $\Pi^*_{[a,b]}(f)$ has a degree less than or equal to two, one would want to compare it with the performance of the 
corresponding Lagrange polynomial with the same degree. \sa
In order to process it, we must assume that $f$ belongs to $C^3([a,b])$, and then, in \cite{Crouzeix_Mignot}, we find that the interpolation error $e_{T}(x)$ for a Lagrange polynomial 
whose degree is less than or equal to two is given by:
\begin{equation}\label{Error_P2}
\forall x \in [a,b]: |e_{T}(x)| \leq \frac{(b-a)^3}{24}M_3,
\end{equation}
where $\D M_3=\sup_{a\leq x\leq b}|f'''(x)|$.\sa
Consequently, by comparing (\ref{Error_P2}) and (\ref{Max2}), provided the given function $f$ is sufficiently smooth, (namely in $C^3([a,b])$), then one would prefer
to use the Lagrange polynomial of degree less than or equal to two, which leads to a more accurate interpolation error. \sa
However, for a function $f$ which only belongs to $C^2([a,b])$, no result is available for this Lagrange polynomial, and the comparison is not valid anymore.
\end{enumerate}
\subsection{The quadrature error}\label{Quadrature_Error}
\noindent We now consider, for any integrable function $f$ defined on $[a,b]$, the famous trapezoidal quadrature \cite{Crouzeix_Mignot}, the formula of which is given by:
\begin{equation}\label{Trapeze}
\D  \int_{a}^{b} f(x)dx \simeq \frac{b-a}{2}\left(f(a)+f(b)\right).
\end{equation}
The reason why we consider (\ref{Trapeze}) is the fact that this quadrature formula corresponds to approximating the function $f$ by its Lagrange polynomial interpolation 
$\Pi_{[a,b]}(f)$, of degree less than or equal to one, which is given by (\ref{P1_Interpolate}). \sa
In the literature on numerical integration, (see for example \cite{Crouzeix_Mignot} and \cite{Cerone}), the following estimation is well known as the \emph{trapezoid inequality} :
\begin{equation}\label{Erreur_Trapeze_Standard}
\D\left|\int_{a}^{b}f(x)\,dx - (b-a)\frac{f(a)+f(b)}{2}\right| \leq \frac{(b-a)^2}{12}\sup_{a\leq x\leq b}|f''(x)|,
\end{equation}
for any function $f$ twice differentiable on $[a,b]$, the second derivative of which is accordingly bounded on $[a,b]$.\sa
It is also well known \cite{Barnett_Dragomir} that if $f$ is only $\mathcal{C}^{1}$ on $[a,b]$, one has the following estimation:
\begin{equation}\label{Erreur_Trapeze_C1}
\D\left|\int_{a}^{b}f(x)\,dx - (b-a)\frac{f(a)+f(b)}{2}\right| \leq \frac{(b-a)^2}{8}(M_1-m_1),
\end{equation}
where $\forall x \in [a,b]: -\infty < m_1 \leq f'(x) \leq M_1 < +\infty$.\sa
So, we will now prove a lemma which will propose estimation (\ref{Erreur_Trapeze_Standard}) in an alternate display. It will also extend estimation (\ref{Erreur_Trapeze_C1}) to twice differentiable functions $f$ which satisfy (\ref{m2M2}).
\begin{lemma}
Let $f$ be a twice differentiable mapping on $[a,b]$ which satisfies (\ref{m2M2}).\sa
Then, we have the following estimation:
\begin{equation}\label{Erreur_Trapeze_C2}
\D\left|\int_{a}^{b}f(x)\,dx - (b-a)\frac{f(a)+f(b)}{2}\right| \leq \frac{(b-a)^3}{24}(M_2-m_2).
\end{equation}
\end{lemma}
\begin{prooff}
In order to derive estimation (\ref{Erreur_Trapeze_C2}), we recall that the classical first-order Taylor formula (\ref{First_Order_Taylor}) enables us to write the 
polynomial $\Pi_{[a,b]}(f)$ by (\ref{Error_P1_Taylor_V0}). Then, by integrating (\ref{Error_P1_Taylor_V0}) between $a$ and $b$, we obtain:
\begin{equation}\label{Erreur_Trapeze_C2_V01}
\D\int_{a}^{b}\left(f(x)\,dx - \Pi_{[a,b]}(f)(x)\right)\,dx = \int_{a}^{b}\left[\frac{(x-a)(b-x)}{b-a}(\epsilon_{x,1}(a)-\epsilon_{x,1}(b))\right]\,dx.
\end{equation}
However, one can easily show that the $P_1-$Lagrange interpolation polynomial $\Pi_{[a,b]}(f)$ given by (\ref{P1_Interpolate}) also fulfills:
\begin{equation}\label{integrale_P1_V0}
\int_{a}^{b}\Pi_{[a,b]}(f)(x)\,dx = \frac{b-a}{2}\left(f(a)+f(b)\right).
\end{equation}
Now, if we introduce the well known quantity $E(f)$, which is called the quadrature error and is defined by:
\begin{equation}\label{erreur_quadrature}
E(f) \equiv \int_{a}^{b} f(x)dx - \frac{(b-a)}{2}(f(a)+f(b)),
\end{equation}
equations (\ref{Erreur_Trapeze_C2_V01}) and (\ref{integrale_P1_V0}) lead to the two following inequalities:
\begin{equation}\label{Mino_Ef}
\D E(f) \geq \frac{1}{2(b-a)}\int_{a}^{b}(x-a)(b-x)\biggl[m_2(x-a)-M_2(b-x)\biggr]\,dx,
\end{equation}
and
\begin{equation}\label{Majo_Ef}
E(f) \leq \frac{1}{2(b-a)}\int_{a}^{b}(x-a)(b-x)\biggl[M_2(x-a)-m_2(b-x)\biggr]\,dx,
\end{equation}
where we used inequality (\ref{Epsilon_b_bounded}) for $\epsilon_{x,1}(a)$ and $\epsilon_{x,1}(b)$, with obvious adaptations.\sa
One can now observe that in (\ref{Mino_Ef}) and (\ref{Majo_Ef}), the two integrals $I$ and $J$ defined by:
\begin{equation}\label{Integral_I_and_J}
\D I = \int_{a}^{b}(b-x)(x-a)^2\, dx \hs \mbox{ and } J = \int_{a}^{b}(x-a)(b-x)^2\, dx
\end{equation}
can be computed as follows. \sa
Let us consider in (\ref{Integral_I_and_J}) the substitution $\D x=\frac{a+b}{2}+\frac{(b-a)}{2}\,t$, then we obtain:
\begin{equation}\label{Integral_I_V1}
\D I = \left(\frac{b-a}{2}\right)^4\int_{-1}^{1}(1-t)(1+t)^2\, dx = \frac{(b-a)^4}{12},
\end{equation}
and
\begin{equation}\label{Integral_J_V1}
\D J = \left(\frac{b-a}{2}\right)^4\int_{-1}^{1}(1+t)(1-t)^2\, dx = \frac{(b-a)^4}{12}.
\end{equation}
Finally, to obtain an upper bound for $|E(f)|$ , owing to (\ref{Mino_Ef})-(\ref{Majo_Ef}) and (\ref{Integral_I_V1})-(\ref{Integral_J_V1}), we obtain:
\begin{equation}\label{Error_Quad_Taylor_Final}
|E(f)| = \left|\int_{a}^{b} f(x)dx - \frac{(b-a)}{2}(f(a)+f(b))\right| \leq  \frac{(b-a)^3}{24}(M_2-m_2).
\end{equation}
\end{prooff}
We now consider the expression of the polynomial interpolation $\Pi_{[a,b]}(f)(x)$ to transform it with the help of the new first order Taylor-like formula (\ref{Taylor_2pts}).
This will enable us to get the following lemma devoted to the corrected trapezoid formula according to Atkison's terminology \cite{Atkinson}.
\begin{lemma}\label{New_Quadrature_Error_Thm}
Let $f$ be a twice differentiable mapping on $[a,b]$ which satisfies (\ref{m2M2}).\sa
Then, we have the following corrected trapezoidal estimation:
\begin{equation}\label{New Erreur_Trapeze_C2}
\D\left|\int_{a}^{b}f(x)\,dx - (b-a)\frac{f(a)+f(b)}{2} + \frac{(b-a)^2(f'(b)-f'a))}{12}\right| \leq \frac{(b-a)^3}{48}(M_2-m_2).
\end{equation}
\end{lemma}
\begin{prooff}
We consider the expression we obtained in (\ref{Error_P1_V00}) for the polynomial interpolation $\Pi_{[a,b]}(f)(x)$, and we integrate it between $a$ and $b$ to get:
\begin{equation}\label{Erreur_Trapeze_C2_V0}
\D\int_{a}^{b}\hspace{-0.3cm}\left[f(x)- \Pi_{[a,b]}(f)(x)\right]dx + \frac{(b-a)^2(f'(b)-f'a))}{12} = \frac{1}{b-a}\int_{a}^{b}(b-x)(x-a)(\epsilon_{x,2}(a)-\epsilon_{x,2}(b))dx,
\end{equation}
where we used the following result obtained by the same substitution which we used to compute integrals in (\ref{Integral_I_and_J}):
$$\D\int_{a}^{b}(b-x)(x-a)\,dx=\frac{(b-a)^3}{6}.$$
Then, due to (\ref{epsilon_x,2(a)_x,2(b)}), we also have the following inequality:
\begin{equation}
\D|\epsilon_{x,2}(a)-\epsilon_{x,2}(b)| \leq \frac{(b-a)}{8}(M_2-m_2),
\end{equation}
and (\ref{Erreur_Trapeze_C2_V0}) directly gives the result (\ref{New Erreur_Trapeze_C2}) to be proved.
\end{prooff}
We conclude this section by several remarks.
\begin{enumerate}
\item We observe that the quadrature error we derived in (\ref{New Erreur_Trapeze_C2}) by the new first-order, Taylor-like formula, is bounded two times less than those we derived in (\ref{Erreur_Trapeze_C2}) by the help of the classical Taylor's formula. \sa
    Furthermore, X.L. Cheng and J. Sun proved in \cite{Xiao} that the best constant one can expect in (\ref{New Erreur_Trapeze_C2}) is equal to $1/36\sqrt{3} \simeq 1/62.35$,
    which is a little bit smaller than $1/48$ we found in (\ref{Erreur_Trapeze_C2_V0}).
\item We notice again that if we consider the particular class of $C^2-$functions $f$ defined on $\R$, $(b-a)$-periodic, so $f'(a)=f'(b),$ and the corrected trapezoid 
formula (\ref{New Erreur_Trapeze_C2}) becomes the classical one:
\begin{equation}
\D\left|\int_{a}^{b}f(x)\,dx - (b-a)\frac{f(a)+f(b)}{2}\right| \leq \frac{(b-a)^3}{48}(M_2-m_2).
\end{equation}
In other words, here we find that for this class of periodic functions, the quadrature error of the classical trapezoid formula is two times more accurate than 
those we found in (\ref{Error_Quad_Taylor_Final}), where the classical first order Taylor formula was implemented.
\end{enumerate}
\section{Conclusions and perspectives}\label{D}
\noindent In this paper we derived a new first order Taylor-like formula to minimize the unknown remainder which appears in the classical one. This new
formula is composed of a linear combination of the first derivative of a given function, computed at $(n+1)$ equally spaced points on $[a,b]$. \sa
We also showed that the corresponding new remainder can be minimized using a suitable choice of the set of the weights which appear in the linear combination
of the first derivative values at the corresponding points.\sa
As a consequence, the new remainder is $2n$ smaller than the one which appears in the classical first order Taylor formula. \sa
Next, we considered two famous applicative contexts given by the numerical analysis where Taylor formula is used: the interpolation error and the quadrature error. Then,
we showed that one can obtain a significative improvement of the corresponding errors. Namely, Lemma \ref{New_error_interpolation_Thm} and Lemma \ref{New_Quadrature_Error_Thm}
prove that the upper bound of these errors is two times smaller than the usual ones estimated by classical Taylor formula, if one limits himself to the class of periodic functions.\sa
Several other applications might be considered by this new first order Taylor-like formula. For example, the approximation error which has to be considered in ODE's where 
Taylor formula is strongly used for appropriate numerical schemes, or also, in the context of finite elements. \sa
For this last application, when one considers linear second elliptic PDE's, due to Cea's lemma \cite{ChaskaPDE}, we noticed that the approximation error is bounded by
the interpolation error. Then, the improvement of the interpolation error we showed in this current work by using the interpolation polynomial
defined by (\ref{corrected_polynomial}), in comparison with the standard $P_1$-Lagrange Polynomial, will consequently impact the accuracy of the approximation error.\sa
Indeed, we highlighted the corresponding gain one may take into account for building meshes, as soon as a given local threshold of accuracy is fixed for the associated approximations.\sa
Other developments may also be considered: It e.g. regards a generalized high order Taylor-like formula, on the one hand, or its corresponding extension for functions with 
several variables, on the other hand.\sa
\textbf{\underline{Homages}:} The authors want to warmly dedicate this research to pay homage to the memory of Professors Andr\'e Avez and G\'erard Tronel, who largely promoted
the passion of research and teaching in mathematics to their students. A special dedication is also expressed to the memory of Mr. Victor Nacasch, who was passionate about the probability theory.
\end{document}